\documentclass[a4paper,12pt]{article}

\usepackage{amsmath}
\usepackage{amssymb}
\usepackage{theorem}

\begin{document}

\begin{center}
INCONSISTENCY OF INACCESSIBILITY
\\
  Alexander Kiselev
\end{center}
\renewcommand{\abstractname}{}
\begin{abstract}
{\footnotesize The work presents the brief exposition of the proof (in \emph{ZF}) of inaccessible
cardinals nonexistence. To this end in view there is used the apparatus of subinaccessible
cardinals and its basic tools ---  reduced formula spectra and matrices and matrix functions and
others. Much attention is devoted to the explicit and substantial development and cultivation of
basic ideas, serving as grounds for all main constructions and reasonings.}
\end{abstract}
\setcounter{equation}{0}
 In 1997 the author has proved

Main theorem ($ZF$):\emph{ There are no weakly inaccessible
cardinals.}

The proof of this theorem was derived as a result of using the
subinaccessible cardinal apparatus which the author has worked out
since 1976, the preliminarily investigations were developed since
1973. The systematic exposition of this proof has been published
in 2000, and the most complete and detailed form the proof of
inaccessible cardinals nonexistence have received in works
``Inaccessibility and Subinaccessibility", Part I~\cite{Kiselev1}
and Part II~\cite{Kiselev2} in 2008, 2010; these two works one can
see also at arXiv sites ~\cite{Kiselev1},~\cite{Kiselev2}.

However, some criticism has been expressed that these works expose
the material which is too complicated  and too extensive and
overloaded by the technical side of the matter, that should be
avoided even when it  uses in essence some new complicated
apparatus . According to these views every result, even  extremely
strong, should be exposed on few pages, otherwise it causes doubts
in its validity.

So, the present work constitutes the brief exposition of the whole
investigation, called to overcome such criticism.

Preliminarily it is required to present the notions of various
subinaccessibles  and reduced matrices and matrix functions  (Part
I and half  of Part II of specified works) - otherwise it is
hardly possible even to sketch the idea of inaccessibles
inconsistency proof.

In these works it is proved that the system
 \quad \   \ $ZF+\exists k \hspace{2mm} (k \; \mbox{\it is weakly inaccessible
    cardinal})$\
\\
is inconsistent. In what follows all the reasoning will be carried
out in this system.

The idea  of main theorem proof consists in the formation of
matrix functions that are sequences of matrices, reduced to a
fixed cardinal; such matrices are certain Boolean values in
L\'{e}vy algebra, reduced to certain cardinal; on this foundations
the simplest matrix function
    \ $S_{\chi f}=( S_{\chi \tau } )_{\tau }$
is introduced  as the sequence of such matrices of special kind
(see below).

This function has domain cofinal to inaccessible cardinal \ $k$,\
it has range consisting of matrices reduced to some fixed cardinal
\ $\chi<k$ \ and defined as minimal in the sense of G\"{o}del
function \ $Od$ \ on corresponding carriers; this property
provides its monotonicity also in the same sense.
 The role of
reducing cardinal \ $\chi$ \ is played further mainly by the
complete cardinal \ $\chi^{\ast}<k$ \ (see below).

Now the idea of the main theorem proof comes out:

\emph{The required contradiction consists in creation of certain matrix function which possess
inconsistent properties: it is monotone and at the same time it is deprived of its monotonicity.}

Let's turn to realization of this idea.

Weakly inaccessible cardinals become strongly inaccessible in
G\"{o}del constructive class \ $L$ \  of values of G\"{o}del
constructive function \ $F$ \ defined on the class of all
ordinals. Every set \ $a\in L$ \ receives its ordinal number \ $Od
( a ) =\min \{ \alpha :F (\alpha  ) =a\}$.
 The starting structure in the further reasoning is
the countable initial segment
    \ $\mathfrak{M}=(L_{\chi ^{0}},\; \in , \; =)$ \quad \
of the class \ $L $ \ serving as the standard model of the theory
    \quad \quad \ \ $ZF+V=L+\exists k \hspace{2mm} (k \; \mbox{\it is
    inaccessible cardinal})$. \quad \

Only the finite part of this theory will be used here because we
shall consider only formulas of limited length (as it will be
clear from what follows). So, the countability of this structure
is required only for some technical convenience and it is possible
to get along without it.

Further \ $k$ \ is the \textit{minimal inaccessible cardinal in} \ $\mathfrak{M}$. We shall
investigate it ``from inside'', considering the hierarchy of subinaccessible cardinals; the latter
are ``inaccessible'' by means of formulas of certain elementary language. To receive this hierarchy
rich enough it is natural to use some rich truth algebra \ $B$. It is well-known that every Boolean
algebra is embedded in an appropriate collapsing \ $ ( \omega _{0},\mu  ) $-algebra and therefore
it is natural to use as \ $B$ \ the sum of the set of such algebras of power \ $k$, \ that is
L\'{e}vy \ $ ( \omega _{0},k ) $-algebra \ $B$. \hfil
\\
This algebra can be introduced in the following way. Let's apply
the set \ $P\in \mathfrak{M}$ \ of forcing conditions that are
finite functions \ $p\subset k\times k$ \ such that for every
limit \ $\alpha <k$ \ and \ $n\in \omega_{0}$ \
    \ $\alpha +n\in dom ( p ) \rightarrow p (\alpha +n ) <\alpha$;
also let \ $p(n)\leq n $ \ for \ $\alpha=0 $. The relation \ $\leq
$ \ of partial order is introduced on \ $P$: \
    \ $p_{1} \leq p_{2} \leftrightarrow p_{2} \subseteq p_{1}$.
After that \ $P$ \ is densely embedded  in the Boolean algebra \
$B\in \mathfrak{M} $ complete in \ $\mathfrak{M}$,\ consisting of
regular sections \ $\subseteq P$.

 It is well known that algebra \ $B$ \ satisfies the \ $k$-chain condition, therefore
 it is possible to consider instead of
values \  $A\in B$ \ only sets
  \ $P_{A}=\{ p\in A:dom (p) \subseteq \chi \}$
where \ $\chi  < k$, \quad \ \ $\chi =\min \{ \chi^{\prime
}:\forall p\in A\quad p\left| \chi ^{\prime }\right. \leq A\}$, \
\quad \ \ (here \ $p\left| \chi ^{\prime }\right.$ \ is the
restriction of \ $p$ \ to \  $\chi ^{\prime }$). Since \  $A=\sum
P_{A}$, \ we shall always identify \  $A$ \ and \  $P_{A}$, \ that
is we shall always consider \ $P_{A}$ \ instead of \ $A$ \ itself.

Just due to this convention and \ $GCH$ in \ $L_k$ \ all Boolean
values \ $A \in B$ \ are sets in \ $L_k$, \ not classes, and this
phenomenon will make possible all further reasoning as a whole.

We shall investigate the hierarchy of subinaccessible cardinals with the help of Boolean values in
\  $B$ \  of some propositions. The countability of the structure \ $\mathfrak{M}$ \  is needed
here only to shorten the reasoning when using its generic extensions by means of \
$\mathfrak{M}$-generic ultrafilters on \ $B$. \ It is possible to get without it developing the
corresponding reasoning in the Boolean-valued universe \ $L^{B}$.

The main instrument of the further reasoning is the notion of a
formula spectrum. Let's use the usual elementary language \
$\mathcal{L}$ \ over the standard structure
    \quad \ \ $( L_{k}\left[ l\right] ,\; \in, \; =,\; l )$,
where \quad \ \emph{l} \quad \ is $\mathfrak{M}$-generic
ultrafilter on \ $B$ .\ Its alphabet consists of usual logic
symbols, all names from Boolean-valued universe \ $L_{k}^{B} $ \
serving as individual constants, symbols \ $\in $, \ $=$  and \
$\underline{l}$,\  the canonical name of\emph{ l}. So, all ordinal
variables  and constants will take values \ $<k$; \ all formulas
will be considered as formulas of the language \ $\mathcal{L}$ \
(if the other case is not meant by the context). Such formulas are
arranged in the elementary Levy hierarchy
  \quad \ \ $\{ \Sigma _{n}( \overrightarrow{a} ) ;\ \Pi _{n}
    ( \overrightarrow{a}) \} _{n\in \omega _{0}}$  \quad \ of formula
classes,
 where  \  $ \overrightarrow{a}$ is a train of individual
constants. Further such classes and their formulas  of some fixed
level \ $n>3$ \  are considered. This agreement is taken to have
in hand further sufficiently large subinaccessible tools and also
to use some auxiliary formulas, terms, relations and sets defined
in \  $L_{k}$ \ directly as additional relational constants in
formulas denotations without raising their level. Obviously, in
this way can be considered \ $P$, \ $B$, relations and operations
on them mentioned above, and also the following:

 \quad $\lessdot $  -- the relation of well-ordering on\quad
$L_{k}$: \ $ \quad a\lessdot b\longleftrightarrow Od ( a ) <Od ( b
) \quad ; $

 \quad $\vartriangleleft $ -- the corresponding relation on \
$L_{k}\times k$:
    \quad \ \ $a \vartriangleleft \beta \longleftrightarrow Od ( a ) < \beta$.

Similarly one can use in \ $L_{k}\left[ l\right] $ \  G\"{o}del
constructive function \ $F^{l}$ \ relatively to \ $l$ \  and
 the function \ $ Od^{l} ( a ) =\min \{
\alpha :F^{l} ( \alpha ) =a\} \ $ and also relations
    \ $a \lessdot ^{l}b\longleftrightarrow Od^{l} ( a ) < Od^{l} ( b
    )~$,
   \ $ a\vartriangleleft ^{l}\beta\longleftrightarrow Od^{l} ( a )
    <\beta~$.

We shall introduce the notion of spectrum only for propositions \
$\varphi
 ( \overrightarrow{a},\underline{l} ) $ \ with train of
individual constants \ $\overrightarrow{a}= ( a_{1},...a_{m} )$ \
consisting of ordinal constants (if the context does not point to
another case). It is possible to manage without this convention
replacing occurrences of each \  $a_{i}$ \ by occurrences of the
term \ $F^{\underline{l}} ( \alpha _{i} ) $ \ for the
corresponding ordinal constant \  $\alpha _{i}$.
\\
Let's also assume that every train \  $\overrightarrow{a}= (
\alpha _{1},...,\alpha _{m} ) $ \  of ordinals \ $ < k$ \ is
identified with the ordinal which is its image under the canonical
order isomorphism of \  $ ^{m}k$ \  onto \  $k$. \

 For every formula \ $\varphi $ \ and
ordinal \ $\alpha _{1}\leq k$ \ by \ $\varphi ^{\vartriangleleft
\alpha _{1}}$  \ is denoted the formula obtained from \  $\varphi$
\ by \ $\vartriangleleft ^{\underline{l}}$~-bounding of all its
quantors by the ordinal \ $\alpha _{1}$, \  that is by replacing
all occurrences of such quantors \  $\exists x$, $\forall x$ \ by
corresponding occurrences of
 \ $\exists x~ ( x\vartriangleleft
^{\underline{l}} \alpha _{1}
    \wedge ... )$, \quad \ \ $\forall x~ ( x\vartriangleleft ^{\underline{l}}
    \alpha _{1}\rightarrow ... )$.
\quad \ In addition, if \  $\alpha _{1}<k$, \  then we say that \
$\varphi $ \  is restricted to \  $\alpha _{1}$ \  or relativized
to \ $\alpha _{1}$; \  if, in addition, the proposition \ $\varphi
^{\vartriangleleft \alpha _{1}}$ \  holds, then we say that \
$\varphi $  \  holds below \  $\alpha _{1}$ \  or that \ $\varphi
$ \  is preserved under this restriction or relativization to \
$\alpha _{1}$. In all such cases  \ $\alpha_{1}$ \ is named
respectively the \ $\vartriangleleft ^{\underline{l}}$~-bounding
ordinal. It is obvious, that for  \  $\alpha _{1}<k$ \ all
formulas \ $\varphi ^{\vartriangleleft \alpha _{1}}$ \ belong to
the class \ $\Delta_{1}$.

If \ $\alpha _{1}=k$, \ then the upper index \ $\vartriangleleft
\alpha _{1}$ \  is omitted and such formulas, reasoning and
constructions are named unrestricted or unrelativized.

\textbf{Definition 1} \quad \ \emph{1) Let \  $\varphi
(\overrightarrow{a},\underline{l} ) $ \  be a proposition \
$\exists x~\varphi _{1} (x,\overrightarrow{a}, \underline{l} ) $
and \ $\alpha _{1}\leq k$. For every \ $\alpha <\alpha _{1}$ \ let
us introduce the following Boolean values:}

\medskip

 \noindent
    \ $A_{\varphi }^{\vartriangleleft \alpha _{1}} ( \alpha ,
    \overrightarrow{a}) = \left\| \exists x\trianglelefteq ^{\underline{l}}
    \alpha ~\varphi_{1}^{\vartriangleleft \alpha _{1}}
    ( x,\overrightarrow{a},\underline{l} ) \right\| $;
    \quad \ \ $\Delta _{\varphi }^{\vartriangleleft \alpha _{1}} ( \alpha ,
    \overrightarrow{a} ) =A_{\varphi }^{\vartriangleleft \alpha _{1}}
    ( \alpha ,\overrightarrow{a} ) - \sum_{\alpha ^{\prime }<\alpha}
    A_{\varphi }^{\vartriangleleft \alpha _{1}} ( \alpha ^{\prime },
    \overrightarrow{a} )$.

\medskip

\emph{2) We name the following function \ $\mathbf{S}_{\varphi
}^{\vartriangleleft \alpha _{1}} ( \overrightarrow{a} ) $ the
spectrum of \  $\varphi $ \  on the point \  $\overrightarrow{a}$
\  below \  $\alpha _{1}$:}

\medskip

   \quad \quad \quad \quad \ \ $\mathbf{S}_{\varphi }^{\vartriangleleft \alpha _{1}}
    ( \overrightarrow{a}) =\{  ( \alpha ,
    \Delta_{\varphi }^{\vartriangleleft \alpha_{1}}
    ( \alpha ,\overrightarrow{a} )  ) :\alpha <\alpha_{1}
    \wedge \Delta _{\varphi }^{\vartriangleleft \alpha _{1}}
    ( \alpha , \overrightarrow{a} ) >0\}$.

\medskip

{\it Projections
    \quad \ \ $dom \left( \mathbf{S}_{\varphi }^{\vartriangleleft \alpha _{1}}
    ( \overrightarrow{a} ) \right), \quad  rng \left (
    \mathbf{S}_{\varphi }^{\vartriangleleft \alpha _{1}}
    ( \overrightarrow{a} ) \right )$
\quad \ are named respectively the ordinal and the Boolean spectra
of \ $\varphi $ \ on the point \ $\overrightarrow{a}$ \ below \
$\alpha _{1}$.

3) If \ $ ( \alpha ,\Delta  ) \in \mathbf{S}_{\varphi
}^{\vartriangleleft \alpha _{1}} ( \overrightarrow{a} ) $, \ then
\  $\alpha $ \  is named the jump ordinal of this formula and
spectra, while \ $\Delta $ \  is named its Boolean value on the
point \  $\overrightarrow{a}$ \ below \ $\alpha _{1}$.

4) The ordinal \ $\alpha _{1}$
 itself is named the carrier of these spectra.}    \hfill $\dashv$

 \noindent
 If a train \  $\overrightarrow{a}$ \  is empty,
then we omit it in notations and omit other mentionings about it.

The investigation of propositions  is natural by means of their
spectra, so one can develop more fine analysis using their
two-dimensional, three-dimensional spectra and so on.

All spectra introduced possess the following simple properties:

\textbf{Lemma 2} \quad \  \emph{Let \ $\varphi $ \ be a
proposition
    \ $\exists x\ \varphi _{1} ( x,\overrightarrow{a},l ),
     \quad \ \varphi_{1}\in \Pi _{n-1}, \quad \
    \alpha _{1}\leq k$, \hfil
\\
then \quad \  \ $\sup dom \left (
\mathbf{S}_{\varphi}^{\vartriangleleft \alpha_{1}} (
\overrightarrow{a} ) \right ) <k $}. \hfill $\dashv$

Here this lemma comes directly from \emph{k}--chain property of
\emph{B}.
 This and all other spectra basic properties will remain at
all their further transformations.

 The so called \textit{universal spectrum} is singled out among all other spectra.  It is well known that
the class \ $\Sigma_{n} ( \overrightarrow{a})$ \ for \  $n>0$ \
contains the formula which is universal for this class; let's
denote it by \ $U_{n}^{\Sigma } ( \mathfrak{n},
\overrightarrow{a},\underline{l} ) $. \  Its universality means
that for any \  $\Sigma _{n} ( \overrightarrow{a} ) $-formula \  $
\varphi  ( \overrightarrow{a},\underline{l} ) $ \  there is a
natural \ $\mathfrak{n}$ \ (the G\"{o}del number of \ $\varphi $\
) such that
    \quad \ \ $\varphi  ( \overrightarrow{a},\underline{l} ) \
    \longleftrightarrow \
    U_{n}^{\Sigma } ( \mathfrak{n},\overrightarrow{a},\underline{l} )$.  The dual formula
universal for $\Pi_{n} ( \overrightarrow{a} ) $  \  is denoted by
$ U_{n}^{\Pi } ( \mathfrak{n},\overrightarrow{a},\underline{ l} )
$. \ For some convenience we shall use \  $U_{n}^{\Sigma }$ \ in
the form \ $\exists x~U_{n-1}^{\Pi } (
\mathfrak{n},x,\overrightarrow{a}, \underline{l} ) $. \ In this
notation the upper indices \ $^{\Sigma} $,~$^{\Pi} $ \ will be
omitted in the case when they can be  restored from the context.

 We name as the spectral universal for the class \
$\Sigma_{n}$ \ formula of level \ $n$ \ the formula \
$u_{n}^{\Sigma } (\overrightarrow{a}, \underline{l} ) $ \ obtained
from the universal formula \ $U_{n}^{\Sigma } (
\mathfrak{n},\overrightarrow{a},\underline{l} ) $ \  by replacing
all occurrences of the variable \ $\mathfrak{n}$ \ by occurrences
of the term \  $\underline{l} ( \omega _{0} ) $.
 The spectral universal for the class $\Pi_{n}$
formula \ $u_{n}^{\Pi } ( \overrightarrow{a},\underline{l} )$ is
introduced in the dual way. Thus we shall use
     \ $u_{n}^{\Sigma } ( \overrightarrow{a}, \underline{l} ) =
     \exists x~u_{n-1}^{\Pi } ( x,\overrightarrow{a},\underline{l}
     )$,
where \ $u_{n-1}^{\Pi} ( x, \overrightarrow{a}, \underline{l} )$ \
is the spectral universal for the class \ $\Pi_{n-1}$ \ formula.

 \quad The values
    \ $A_{\varphi}^{\vartriangleleft \alpha _{1}} ( \alpha ,
    \overrightarrow{a} ), \Delta _{\varphi }^{\vartriangleleft
    \alpha _{1}} ( \alpha,\overrightarrow{a} )$ and the spectrum
     \ $ \mathbf{S}_{\varphi}^{\vartriangleleft
    \alpha _{1}} ( \overrightarrow{a} )$
\quad \ of the formula \ $\varphi =u_{n}^{\Sigma }
(\overrightarrow{a}, \underline{l} )$ \ and its projections will
be named the universal Boolean values and spectra of the level \
$n$ \  on the point \ $\overrightarrow{a}$ \  below \ $\alpha
_{1}$ \ and they are denoted by

  \quad \quad \quad \quad \quad \  \ $A_{n}^{\vartriangleleft \alpha _{1}} ( \alpha ,
    \overrightarrow{a} ), \quad \Delta _{n}^{\vartriangleleft
    \alpha_{1}} ( \alpha ,\overrightarrow{a} )$,\quad
    \ $\mathbf{S}_{n}^{\vartriangleleft \alpha _{1}}
    ( \overrightarrow{a} )$.

 \noindent
Here the term ``universal spectra'' is justified by the fact:
 for every \  $\varphi =\exists x\
\varphi _{1} ( x, \overrightarrow{a},\underline{l} )$ ,\quad \
$\varphi _{1}\in \Pi _{n-1}$ there holds
   \quad \ \ $ dom \left (
    \mathbf{S}_{\varphi }^{\vartriangleleft \alpha _{1}} (
    \overrightarrow{a} )\right ) \subseteq dom \left (
    \mathbf{S}_{n}^{\vartriangleleft \alpha _{1}}
    (\overrightarrow{a} ) \right)$ .

All further reasoning is conducted in \  $L_{k}$ \  (or in  \
$\mathfrak{M} $ \ ). Let's introduce the central notion of
subinaccessibility -- the inaccessibility by means of our
language. The ``meaning'' of propositions is contained in their
spectra and therefore it is natural to define this inaccessibility
by means of the spectra of all propositions of a given level:

\textbf{Definition 3} \quad \emph{\ Let \  $\alpha _{1} \leq k $.
We name an ordinal  \ $\alpha <\alpha _{1}$ \  subinaccessible of
\ $a $ \ level \ $n$ \ below \  $\alpha _{1}$ \  iff it fulfills
the
 formula
\quad \ \ $\forall \overrightarrow{a} <\alpha \quad dom
    \left ( \mathbf{S}_{n}^{\vartriangleleft \alpha_{1}}
    ( \overrightarrow{a} ) \right ) \subseteq \alpha $
\quad \ denoted by \  $SIN_{n}^{<\alpha _{1}} ( \alpha)$. The set
    \quad \ \ $\{ \alpha <\alpha _{1}:SIN_{n}^{<\alpha _{1}} ( \alpha )
    \}$
 \quad \ of all these ordinals is denoted by \ $SIN_{n}^{<\alpha _{1}}$ \
and they are named \ $SIN_{n}^{<\alpha _{1}}$-ordinals.}\hfill
$\dashv$

As usual, for \ $\alpha _{1}<k$ \  we say that subinaccessibility
of \  $\alpha $ \  is restricted by \  $\alpha _{1}$ \ or
relativized to \  $\alpha _{1}$; \ for \ $\alpha_1 = k$ \ the
upper indices \ $< \alpha_1$, $\vartriangleleft \alpha_1$ \ are
dropped. This definition obviously equivalent to the following:

 \noindent
 let \  $\alpha <\alpha
_{1}\leq k$, \ $\alpha \in SIN_{n}^{<\alpha _{1}}$  \ and a
proposition \ $ \exists x~\varphi  ( x,\overrightarrow{a},
\underline{l} ) $ \ has \  $ \overrightarrow{a}<\alpha $, \
$\varphi \in \Pi _{n-1}$, \ then for any \ $\mathfrak{M}$-generic
\  $l$
 \quad \ \ $L_{k}\left[ l\right] \vDash
\left( \exists x\vartriangleleft ^{l}
    \alpha _{1}\mathit{\ }\varphi ^{\vartriangleleft \alpha _{1}}
    ( x,\overrightarrow{a}, l ) \longrightarrow \exists x
    \vartriangleleft ^{l}\alpha \ \varphi ^{\vartriangleleft
    \alpha _{1}} ( x,\overrightarrow{a},l) \right)$.\

 \noindent
 In this case we say that below \ $ \alpha _{1}$ \  the ordinal
 \ $\alpha $ \  restricts or relativizes the proposition \ $\exists
x~\varphi $.  Considering the same in the inverted form for \
$\varphi \in \Sigma _{n-1}$:
 \ $L_{k}\left[ l\right] \vDash
( \forall x\vartriangleleft ^{l}
    \alpha \ \varphi ^{\vartriangleleft \alpha _{1}} ( x,
    \overrightarrow{a},l ) \longrightarrow$ \ $ \longrightarrow \forall x
    \vartriangleleft ^{l}\alpha _{1}\varphi ^{\vartriangleleft
    \alpha _{1}} ( x,\overrightarrow{a},l ))$,\
we say that below \  $\alpha _{1}$ \  the ordinal \  $\alpha $ \
extends or prolongs the proposition \  $\forall x\ \varphi $ \ up
to \ $\alpha _{1}$.

Obviously, the cardinal \ $k$ \ is \textit{subinaccessible} itself
of any level, if we define this notion for \ $\alpha = \alpha_1 =
k$.
 So, the comparison of inaccessibility and
subinaccessibility notions naturally arises in a following way:

 \noindent
 The cardinal \ $k$ \ is weakly inaccessible, since it is
uncountable and cannot be reached by means of smaller powers in
the sense that: 1) it is regular and 2) it is closed under
operation of passing to next power: \ $\forall \alpha < k ~~~
\alpha^+ < k$.
  Turning to \textit{subinaccessibility} of the ordinal \ $\alpha < k$ \ of the
level \ $n$, \ one can see that the property of regularity is
dropped now, but \ $\alpha$ \ still can not be reached, but by
more mighty means: condition 2) is strengthened  and \ $\alpha$ \
is closed under more mighty operations of passing to jump ordinals
of universal spectrum:
    \quad \ \ $\forall \overrightarrow{a} <\alpha$ \quad \ $\forall \alpha^{\prime}
    \in dom \mathbf{S}_{n} ( \overrightarrow{a} )
    \quad \alpha^{\prime} < \alpha$,
that is by means of ordinal spectra \textit{of all} propositions
of level \ $n$. \ Hence, this ordinal \ $\alpha$ \ is closed
\textit{under all \ $\Pi_{n-1}$-functions in all extensions} \quad
\ \ $( L_{k}\left[ l\right] ,\; \in, \; =,\; l )$, not only under
operation of power successor in \ $L_k$ .\ In particular, for
every \  $n \geq 2$
    \ $ \alpha =\omega _{\alpha } \mathrm{\quad (in\ } L_{k} )$.
Besides that the set
 \  $SIN_{n}^{<\alpha _{1}}$ \  is closed in \
$\alpha _{1}$, \ that is for any \  $\alpha <\alpha _{1}$
  \quad \  \ $\sup  ( \alpha \cap SIN_{n}^{<\alpha _{1}} ) \in
    SIN_{n}^{<\alpha_{1}}$,
 and the set \  $SIN_{n}$ \  is unbounded in \ $k$, \ $\sup
SIN_{n}=k$.

When formulas are equivalently transformed their spectra can
change. It is possible to use this phenomenon for the analysis of
subinaccessible cardinals. To this end we shall introduce the
universal formulas with ordinal spectra containing only
subinaccessible cardinals of smaller level. For more clearness
formulas without individual constants will be considered. Let's
start with the spectral universal formula for the class \ $\Sigma
_{n}$. \  The upper indices \ $^{\Sigma} $, $ ^{\Pi}$ \ will be
omitted as usual (if it will
not cause misunderstanding). \\
In what follows bounding ordinals \ $\alpha $ \ are always assumed
to be \ $SIN_{n-2}$-ordinals or \ $\alpha =k$.\

\textbf{ Definition 4 } \emph{1) We name as the monotone spectral
universal for the class \  $\Sigma _{n}$ \ formula of the level \
$n$ \  the \ $\Sigma _{n}$-formula
    \ $\widetilde{u}_{n} ( \underline{l} ) =\exists x~\widetilde{u}_{n-1}
    ( x, \underline{l} )$
  \ where \  $\widetilde{u}_{n-1} ( \underline{l} ) \in \Pi _{n-1}$
\  and}
 \ $\widetilde{u}_{n-1} ( x,\underline{l} ) \longleftrightarrow
    \exists x^{\prime}\vartriangleleft
    ^{\underline{l}}x~u_{n-1}^{\Pi}
    ( x^{\prime },\underline{l} )$ .

2) \quad\emph{ We name as the subinaccessibly universal for the
class \ $\Sigma _{n}$ \  formula of the level \  $n$ \  the \
$\Sigma _{n}$-formula
    \quad \ \ $\widetilde{u}_{n}^{\sin } ( \underline{l} ) =
    \exists x~\widetilde{u}_{n-1}^{\sin} ( x,\underline{l} )$
\ where \ $\widetilde{u}_{n-1}^{\sin }\in \Pi _{n-1} $  \  and}
 \ $\widetilde{u}_{n-1}^{\sin } ( x,\underline{l} )
    \longleftrightarrow$ \ $ \longleftrightarrow SIN_{n-1} ( x )
    \wedge \widetilde{u}_{n-1} ( x ,
    \underline{l} ) $.

\emph{The subinaccessibly universal for the class \ $\Pi_{n} $ \
formula is introduced in the dual way.}

\emph{3) \quad The Boolean values
    \quad \ \ $A_{\varphi}^{\vartriangleleft \alpha _{1}} ( \alpha  ), \
    \Delta_{\varphi }^{\vartriangleleft \alpha _{1}} ( \alpha  ),
    \mathrm{\it \ \emph{\emph{ and the spectrum  }}}
    \mathbf{S}_{\varphi }^{\vartriangleleft \alpha_{1}}$
of the formula \ $\varphi = \widetilde{u}_{n}^{\sin }$ \ and its
projections (see definition 1 where \
$\trianglelefteq^{\underline{l}}$ \ should be replaced with \
$\leq$) are named subinaccessibly universal of the level \ $n$ \
below \ $\alpha_{1}$ \ and are denoted respectively by
  \quad \ \ $\widetilde{A}_{n}^{\sin \vartriangleleft \alpha _{1}} ( \alpha ) ,
    \quad \widetilde{\Delta }_{n}^{\sin \vartriangleleft \alpha _{1}}
    ( \alpha  ) ,\quad \widetilde{\mathbf{S}}_{n}^{\sin \vartriangleleft
    \alpha _{1}}$.}        \hspace*{\fill} $\dashv$

 \noindent
Obviously \quad \ \ $u_{n}^{\Sigma} ( \underline{l} )
\longleftrightarrow \ \widetilde{u}_{n}^{\sin } ( \underline{l}
)$~ \quad \ and
 \quad \ \ $dom \bigl ( \widetilde{\mathbf{S}}_{n}^{\sin \vartriangleleft \alpha
_{1}} \bigr ) \subseteq SIN_{n-1}^{<\alpha_{1}} \cap dom  \bigl
(\mathbf{S}_{n}^{\triangleleft \alpha_{1}} \bigr )$~.

Our aim is to ``compare'' universal spectra with each other on
\textit{different} carriers \ $\alpha_1$ \ disposed cofinally to \
$k$ \ in order to introduce \textit{monotone} matrix functions. To
this end it is natural to do it by means of using values of
function \ $Od$ \ for such spectra.  Also it is natural to find
some estimates of ``informational complexity'' of these spectra by
means of estimates of their order types.
 But the required comparison of such spectra can be hardly carried out in a
proper natural way because they \textit{are ``too much differ''}
from each other for arbitrary great carriers \ $\alpha_1$.

 \noindent
So, there is nothing for it but to consider further
\textit{spectra reduced to some fixed cardinal} and, next,
\textit{reduced matrices}. So here we start to form the main
material for building matrix functions -- reduced matrices. With
this end in view first we shall consider the necessary preliminary
constructions -- reduced spectra.

For an ordinal \ $\chi \leq k$ \ let \ $P_{\chi }$ \ denote the
set \  $\{ p\in P:dom ( p ) \subseteq \chi \} $ \ and \  $B_{\chi
}$ \  denote the subalgebra of \  $B$ \  generated by \ $P_{\chi}$
\  in \  $L_{k}$. \  For every \  $A\in B$ \  let's introduce the
set
    \quad \ \ $A \lceil \chi =\{ p\in P_{\chi }:\exists q~ ( p=\left. q\right|
    \chi \wedge q\leq A ) \}$
 \quad \ which is named the value of \ $A$ \  reduced to \  $\chi $. \ It
is known that
    \quad \ \ $B_{\chi }=\left \{ \sum X:X\subseteq P_{\chi } \right
    \}$
 \quad \ and therefore every \ $A\in B_{\chi }$  \ coincides with \ $ \sum
A\lceil \chi $. \ Therefore let's identify every \ $A\in B_{\chi
}$ \ with its reduced value \ $A\lceil \chi $; \ so, here one
should point out again, that due to the chain property of \ $B$
and \ $GCH$ in  \ $L_k$  every value \ $A \in B_\chi$ \ is the set
in \ $L_k$, not class, and \ $B_\chi$ \ is considered as the set
of such values, \ $ B_\chi \in L_k$
 for \ $\chi <k$. \

\textbf{Definition 5}\quad \  \emph{1) For every \  $\alpha
<\alpha _{1}$  \  let us introduce the Boolean values and the
spectrum:}

    \quad \quad \quad \quad \ \ $\mathbf{\ }\widetilde{A}_{n}^{\sin \vartriangleleft \alpha _{1}}
    (\alpha  ) \lceil \chi \mathbf{;~} \quad \ \widetilde{\Delta }_{n}^{\sin
    \vartriangleleft \alpha _{1}} ( \alpha  ) \overline{\lceil}\chi =
    \widetilde{A}_{n}^{\sin \vartriangleleft \alpha _{1}} ( \alpha  )
    \lceil\chi \mathbf{-}\sum_{\alpha ^{\prime }<\alpha }
    \widetilde{A}_{n}^{\sin \vartriangleleft \alpha _{1}}
    ( \alpha ^{\prime } ) \lceil \chi $;

\[
    \widetilde{\mathbf{S}}_{n}^{\sin \vartriangleleft \alpha _{1}}
    \overline{\overline{\lceil}}\chi =\{
    ( \alpha ;~\widetilde{\Delta }_{n}^{\sin \vartriangleleft \alpha _{1}}
    ( \alpha  ) \overline{\lceil}\chi  ) :\alpha <\alpha _{1}\wedge
    \widetilde{\Delta }_{n}^{\sin \vartriangleleft \alpha _{1}} ( \alpha  )
    \overline{\lceil}\chi >0\}  .
\]
2)\quad \emph{ These values, spectrum and its  first and second
projections are named subinaccessibly universal reduced to} \
$\chi $ \  \emph{of the level \ $n$\quad below} \  $\alpha _{1}$.\
\hspace*{\fill} $\dashv$

 \noindent
 In a similar way multi-dimensional reduced spectra can be introduced.

 \noindent
 From this definition comes that reduced spectra possess previous properties, for
instance, for
 \ $\alpha <\alpha _{1}$,\quad $ \chi \leq k$ \ there holds
  \quad $ \sup dom  \left (
\widetilde{\mathbf{S}}_{n}^{\sin \vartriangleleft
\alpha_{1}}\overline{\overline{\lceil }}\chi \right) < k$, \
 \ $dom \left ( \widetilde{\mathbf{S}}_{n}^{\sin \vartriangleleft
\alpha_{1}}\overline{\overline{\lceil }}\chi \right)
 \subseteq
SIN_{n-1}^{<\alpha_{1}} \cap$ \ $ \cap dom \bigl
(\mathbf{S}_{n}^{\triangleleft \alpha_{1}} \bigr )$~ and so on.
\\
The main role further is played by matrices and spectra reduced to
complete cardinals; their existence comes out from \emph{k}-chain
property of algebra B:

 \textbf{Definition  6} \quad \ \emph{We name as a complete
ordinal of level \ $n$ \  every ordinal \  $\chi $ \  such that}

 \noindent
 \ $\exists x~\widetilde{u}_{n-1}^{\sin } ( x,\underline{l} )
    \longleftrightarrow \ \exists x<\chi
    ~\widetilde{u}_{n-1}^{\sin}( x,\underline{l} )$ .
\quad \ \emph{The least of these  ordinals is denoted by} \  $\chi
^{\ast }$. \ \hspace*{\fill} $\dashv$

 \noindent
From this definition it comes
 \ $\chi ^{\ast }=\sup
dom \left ( \widetilde{\mathbf{S}}_{n}^{\sin } \right ) =\sup dom
\left ( \mathbf{S}_{n}\right ) <k$ \ and \ $SIN_{n-1} ( \chi
^{\ast } )$  and so on.

Let's turn to order spectrum types. If \ $X$ \ is a well ordered
set, then its order type is denoted by \ $OT ( X ) $; \  if \  $X$
\  is a function with well ordered domain, then we assume \  $OT (
X ) =OT ( dom ( X )  ) $. The obvious rough upper estimate of
spectrum types comes from  \ $\left| P_{\chi }\right| =\left| \chi
\right| $ \ and \  $GCH$ \ in \ $L_{k}$~:
    \ $OT ( \widetilde{\mathbf{S}}_{n}^{\sin \vartriangleleft
    \alpha_{1}}\overline{\overline{\lceil }}\chi  ) < \chi
    ^{+}$.\

 \noindent
 \quad \ Now let's discuss estimates of such types from below. Here comes
out the lemma essential for the proof of main theorem: it shows, that as soon as an ordinal \ $
\delta < \chi ^{\ast +}$ \ is defined through some jump ordinals of the subinaccessibly universal
spectrum reduced to \  $\chi ^{\ast }$, the order type of this spectrum exceeds  \ $\delta$ \ under
certain natural conditions:

\textbf{Lemma 7 }{\em (About spectrum type)}  \emph{Let ordinals \
$ \delta$,~$\ \alpha_{0}$,~$\ \alpha_1$ \ be such that:}

 \noindent
(i) \quad $\ \delta < \chi ^{\ast+} <\ \alpha_{0}<\ \alpha_1\leq
k$~;
 \quad \quad (ii)  $SIN_{n-2} ( \ \alpha_1 ) \wedge OT (
SIN_{n-1}^{<\ \alpha_1 } ) =\ \alpha_1$~;

 \noindent
 (iii)\ $\ \alpha _{0}\in dom  \left (
\widetilde{\mathbf{S}} _{n}^{\sin \vartriangleleft \ \alpha_1}\overline{\overline{\lceil }} \
\chi^{\ast} \right ) $; (iv) \ $ \delta $ \ {\it is defined in  \ $ L_{k}$ \ through \ $ \alpha
_{0}$, \ $\chi ^{\ast }$ by a formula \ $\in \Sigma _{n-2}\cup \Pi _{n-2}$~.

 \noindent
 Then} $\quad \quad \quad \quad \quad \quad
\quad \quad \quad \quad\ \delta <OT (
\widetilde{\mathbf{S}}_{n}^{\sin
\vartriangleleft \ \alpha_1}\overline{\overline{\lceil }}%
\ \chi^{\ast }  ) \quad.$          \hspace*{\fill} $\dashv$

 But
still there is the following essential inconvenience: such
spectra, taken on their different carriers, can be hardly compared
with each other in view to their basic properties, because their
domains (ordinal spectrums) can contain an arbitrary great
cardinals, when these carriers are increasing up to \ $k$. In
order to avoid this obstacle we shall transform them to reduced
matrices. These matrices comes from reduced spectra by easy
isomorphic enumeration of their domains:

\textbf{Definition 8} \quad {\it 1) We name as a matrix reduced to an ordinal \ $\chi $ \  every
function \  $M$ \  defined on some ordinal and with  \ $ \ rng ( M ) \subseteq B_{\chi}$ .}

 \noindent
2) \quad \emph{Let \  $M$ \  be a matrix and \ $M_{1}\subset
k\times B$. We name as a superimposition of \ $M$ \ onto \ $M_{1}$
\    an order isomorphism \quad \  f  \quad \ of \ $\ dom (M )$ \
onto \ \ $\ dom ( M_{1} )$ such that}

\quad \quad \quad \  \ $ \forall \alpha ,\alpha ^{\prime }\forall
\Delta ,\Delta ^{\prime } ( f ( \alpha ) =\alpha ^{\prime }\wedge
( \alpha ,\Delta  ) \in M\wedge  ( \alpha ^{\prime },\Delta
^{\prime } ) \in M_{1}\longrightarrow \Delta =\Delta ^{\prime }
)$;

 \noindent
 {\it in this case we say that \  $M$
 \  is superimposed onto \  $M_{1}$ \  and write \
$M\Rightarrow M_{1} $.}

 \noindent
 {\it 3) \quad Let matrix \  $M$ \ be superimposed onto the spectrum \ $\widetilde{\mathbf{S}}_{n}^{\sin
\vartriangleleft \alpha }\overline{\overline{\lceil }}\chi $ then \ $M$ \  is named the matrix of
this spectrum or the subinaccessibly universal matrix of the level \  $n$ \ reduced to \ $\chi$ \
on \ $\alpha $.

 \noindent
4)\quad \ If \ $(\alpha^{\prime} ,\Delta ) \in
\widetilde{\mathbf{S}}_{n}^{\sin \vartriangleleft \alpha
}\overline{\overline{\lceil }}\chi $, \ then  \ $\alpha^{\prime}$
\ is named the jump cardinal of the matrix \ $M$, \ while \
$\Delta$ \ is named its Boolean value on \ $\alpha$.}\quad \
\emph{5) In this case the cardinal \ $\alpha$ \ is named the
carrier of the matrix \ $M$.}            \hspace*{\fill} $\dashv$

It is obvious that for
 \ $M\Rightarrow
\widetilde{\mathbf{S}}_{n}^{\sin
    \vartriangleleft \alpha } \overline{\overline{\lceil }}\chi$
 \quad \ there holds  \ $OT ( M ) =dom ( M ) \leq Od ( M ) <\chi
^{+}.$

 \noindent
 This property  shows that reduced matrices can
be compared (in the sense of function $Od$) \textit{within \
$L_{\chi^{+}}$ \ only} and it will help to define matrix functions
with required properties.

 \noindent
Further the very special role is played by the so called
\textit{singular matrices}:

\textbf{ Definition 9 } \quad \emph{\ We denote by \  $\sigma (
\chi ,\alpha ) $ \ the following formula:}

 \emph{\ $SIN_{n-2} ( \alpha  ) \wedge  ( \chi\mbox{ is a limit cardinal
}<\alpha ) \wedge \ OT ( SIN_{n-1}^{<\alpha } ) =\alpha \wedge
\sup dom \left ( \widetilde{\mathbf{S}}_{n}^{\sin \vartriangleleft
\alpha }\overline{\overline{\lceil}}\chi\right ) =\alpha $~.}

 \noindent
\emph{And let \  $\sigma  ( \chi ,\alpha ,M ) $ \  denote the
formula
    \ $\sigma  ( \chi ,\alpha  ) \wedge  ( M\Rightarrow
    \widetilde{\mathbf{S}}_{n}^{\sin \vartriangleleft \alpha }
    \overline{\overline{\lceil }}\chi  )$ .
The matrix \  $M$  \  and the spectrum \
$\widetilde{\mathbf{S}}_{n}^{\sin \vartriangleleft \alpha}
\overline{\overline{\lceil }}\chi $ \  reduced to \  $\chi $ \ are
named singular on a carrier \ $\alpha $ \   \textit{iff} \ $\sigma
 ( \chi ,\alpha ,M ) $ \  is fulfilled .
 The symbol \ $S$  \  will be used for the common denotation of singular
matrices. }            \hspace*{\fill} $\dashv$

 \noindent
Further all matrices will be singular on their carriers; all
reasoning will be conducted in \ $L_{k}$\ (or in \
$\mathfrak{M}$).

 \noindent
 By this definition jump cardinals
of singular matrix on its carrier \ $\alpha $ \ are disposed
cofinally to \ $\alpha $. \ Due to the last fact it is possible to
introduce the following important cardinals:

\textbf{Definition 10 }   \quad \ \emph{Let \ $\sigma  ( \chi
,\alpha, S ) $ \ fulfills, then we name as jump cardinal and
prejump cardinal after \ $\chi $ \ of the matrix \ $S$ \ on the
carrier \ $\alpha $, \ or, briefly, of the cardinal \ $\alpha$, \
the following cardinals respectively:}

\quad \quad \quad \quad \ \
$\alpha _{\chi }^{\downarrow }=\min \{ \alpha ^{\prime}
    \in ]\chi ,\alpha \lbrack ~:~ \widetilde{A}_{n}^{\sin
    \vartriangleleft \alpha } ( \chi  ) \lceil \chi
    <\widetilde{A}_{n}^{\sin \vartriangleleft \alpha } (
    \alpha^{\prime } ) \lceil \chi \wedge SIN_{n-1}^{<\alpha } ( \alpha^{\prime } ) \}$ ;

\medskip

\quad \quad \quad \quad \ \ $\alpha _{\chi }^{\Downarrow }=\sup \{
\alpha ^{\prime}
    <\alpha _{\chi }^{\downarrow }  ~:~
    \widetilde{A}_{n}^{\sin \vartriangleleft \alpha } ( \chi )
    \lceil \chi =\widetilde{A}_{n}^{\sin \vartriangleleft \alpha }
    (\alpha ^{\prime }) \lceil \chi \wedge SIN_{n-1}^{<\alpha } ( \alpha ^{\prime } ) \}$.\  \hspace*{\fill} $\dashv$

 \noindent
 It is not hard to see, that in this situation these
 \ $\alpha _{\chi }^{\downarrow }$,\  \ $\alpha _{\chi }^{\Downarrow }$\ really do exist  and

 \noindent
 \ $~\alpha _{\chi }^{\downarrow }=\min \left\{ \alpha
^{\prime }>\chi :\alpha ^{\prime }\in dom   \left (
\widetilde{\mathbf{S}}_{n}^{\sin \vartriangleleft \alpha }
\overline{\overline{\lceil }}\chi \right ) \right\}$; \ \quad \ \
$\alpha _{\chi }^{\Downarrow }<\alpha _{\chi }^{\downarrow
}<\alpha ;\quad \quad ]\alpha _{\chi }^{\Downarrow },\alpha _{\chi
}^{\downarrow }[\cap SIN_{n-1}^{<\alpha }=\varnothing$. \

The basic instruments of main theorem proof are matrix functions
that are sequences of reduced singular matrices. The following
lemma makes it possible to build such functions:

\textbf{Lemma 11}
    \quad \ \ $\forall \chi \forall \alpha _{0} (  ( \chi
    \emph{\mbox{\ is a limit cardinal }}>\omega _{0} ) \longrightarrow
    \exists \alpha _{1}>\alpha _{0}~\sigma  ( \chi ,\alpha _{1} )  )$
    . \hspace*{\fill} $\dashv$

 The building of matrix functions relies on the following
enumeration (in~~$L_{k}$) \  of subinaccessible
 \ $SIN_{n-1}$-cardinals below $ \alpha _{1}$:

\textbf{Definition 12} \emph{ The function \ $\gamma
_{f}^{<\alpha_{1}}= ( \gamma _{\tau }^{<\alpha _{1}} ) _{\tau }$~
is defined by recursion on} \  $ \tau <\alpha _{1}\leq k $:

   \ $\gamma _{0}^{<\alpha _{1}}=0 $;  \quad \ \emph{for} \ $\tau>0 $\
 \quad \ \ $ \gamma _{\tau }^{<\alpha
_{1}}=\min \{ \gamma <\alpha_{1}:
    SIN_{n-1}^{<\alpha _{1}} ( \gamma  ) \wedge \forall
    \tau^{\prime }<\tau \quad
    \gamma_{\tau ^{\prime }}^{<\alpha_{1}}<\gamma \}$. \         \hspace*{\fill} $\dashv$

 \noindent
 Obviously, if \ $ \alpha _{1}=k$ then range and domain of this function are both
cofinal to \ $k$.

 The proof of main theorem consists in creation
in \ $L_{k}$ \ of the special matrix functions possessing
inconsistent properties; here are such functions of the main kind:

\textbf{Definition 13 } {\it \ We name as a matrix function of the level \ $n$ \ below $ \alpha
_{1}$ reduced to $\chi $ the function
 \ $S_{\chi f}^{<\alpha _{1}}= ( S_{\chi \tau
}^{<\alpha _{1}} ) _{\tau }$ \

 \noindent
taking values for \ $\tau<k$: \quad \quad \ \ $ S_{\chi \tau
}^{<\alpha _{1}}=\min_{\underline{\lessdot }}\{
    S:\exists \alpha < \alpha_{1}  \left ( \gamma_{\tau}^{<
    \alpha_{1}}< \alpha \wedge \sigma ^{\vartriangleleft
    \alpha _{1}} ( \chi ,\alpha , S )\right ) \}$.}       \hspace*{\fill} $\dashv$

 \noindent
 So, these values are matrices  \ $S_{\chi \tau }^{<\alpha _{1}}$ reduced to \
$\chi$ \ singular on these carriers \ $\alpha$ that are
$\underline{\lessdot }$-minimal for \ $ \gamma_{\tau}^{<
    \alpha_{1}}< \alpha$.  In the same sense all these
values are bounded by  \ $\chi ^{+}<k$ due to \emph{GCH} in \
$L_k$ : \ $ Od (S_{\chi \tau }^{<\alpha _{1}}) <\chi ^{+}.$

Everywhere further \  $\chi =\chi^{\ast} $; the lower index \
$\chi^{\ast} $ \ can be omitted in notations, for instance \
$\alpha _{\chi^{\ast}}^{\downarrow }$,\ \ $\alpha
_{\chi^{\ast}}^{\Downarrow }$,\ \ $S_{\chi^{\ast} f}$, \ \
$S_{\chi ^{\ast }\tau }$ will be denoted through \
$\alpha^{\downarrow }$,\ \ $\alpha^{\Downarrow }$,\ \ $S_{ f}$, \
$S_{\tau }$ (if the context will not mean another case) and so on.

The following two lemmas represent the mainstream of all further
reasoning: they establish that such matrix functions have
properties of definiteness and of $\underline{\lessdot
}$-monotonicity which comes from $\underline{\lessdot
}$--minimality of their values. From lemma 11 there follows
directly:

\textbf{Lemma 14} \emph{(About matrix function definiteness)}
\emph{ The unrelativized function \ $S_{\chi^{\ast} f}$ \ is
defined on the final segment of \ $k$:}
    \quad \ \ $dom \left( S_{\chi ^{\ast } f } \right)=
    \left \{ \tau : \chi^{\ast} \leq \gamma_{\tau} < k \right \}$.   \hspace*{\fill} $\dashv$

 \textbf{Lemma 15}  \emph{This function  \ $S_{\chi^{\ast} f}$ \  is $\underline{\lessdot
}$-monotone:}
    \ $\quad \chi^{\ast} < \tau _{1}<\tau _{2} < k \longrightarrow S_{\chi^{\ast}
    \tau _{1}} \underline{\lessdot }
    S_{\chi^{\ast} \tau _{2}}$. \ \hspace*{\fill} $\dashv$

  Hence \quad this function stabilizes, that is there is an ordinal
\ $\tau ^{\ast }$ \ such that for every \ $\tau \geq \tau ^{\ast
}$ there exist constant value
    \quad \ \ $S_{\chi ^{\ast }\tau }=S_{\chi ^{\ast }\tau ^{\ast }}$.
 \quad \ Thus for all \ $\gamma_{\tau} \geq \chi ^{\ast }$
 \
values  \ $S_{\chi^{\ast} \tau}$    are bounded by the fixed
ordinal \ $Od \left ( S_{\chi ^{\ast }\tau^{\ast}} \right ) <
\chi^{\ast +} $: \
 \quad \quad \   \ $S_{\chi ^{\ast }\tau } \ \underline{\lessdot}
    \ S_{\chi ^{\ast }\tau ^{\ast }} \vartriangleleft \chi^{\ast +}$.

\medskip

Now everything is ready to present the main theorem proof mode.
The monotonicity of the simplest matrix function \ $S_{\chi^{\ast}
f}$ \ is established already. So, the required contradiction
should be find out in its {\sl nonmonotonicity} on the following
way:

 \noindent
The lower index \ $\chi ^{\ast } $ \ will be dropped for some
brevity.
 Let's consider the matrix function  \ $S_f$ in its state of
 stabilizing, that is consider \ an arbitrary sufficiently great
\
 \ $\tau _{0}>\tau ^{\ast }$ \ and the matrix \  $ S_{\tau _{0}}$ \
\ on some carrier \ $\alpha _{0}\in \left] \gamma_{\tau _{0}},k
 \right[ $,
\ the prejump cardinal \ $\alpha^0=\alpha _{0}^{\Downarrow }$ \
and the relativized function \  $ S_{f}^{<\alpha^0}$  below \
$\alpha^0$. \ Remind, all values of \ $S_f$  are bounded by the
fixed ordinal \ $Od \left ( S_{\tau^{\ast}} \right ) < \chi^{\ast
+} $. \

Standing on  the jump cardinal \ $\alpha _{0}^{\downarrow }$
 after \ $\chi^{\ast}$ \ of this matrix on this carrier, one should observe
the behavior of this very function, but in its relativized to \
 \ $\alpha^0=\alpha _{0}^{\Downarrow }$ \ form \
$S_{f}^{<\alpha^0}$ \ below \ $\alpha^0$. \ The function \ $S_f$
is monotone and this relativized function \ $S_{f}^{<\alpha^0}$ \
is monotone too by the same reasons.

But it is excluded. To see it one should apply lemma 7 about
spectrum type to this situation, considering it for
    \quad \ \ $\ \delta = \mathrm{sup}_{\tau} Od
    (S_{\tau }^{<\alpha^0} )$, \quad \ $\alpha_{0}$ - the jump
    cardinal
    \ $\alpha _{0}^{\downarrow }$, \quad \ \ $\alpha_1$ -- the carrier
    \ $\alpha_{0}$.

 \noindent
Let's consider \ $ \delta <\chi ^{\ast +}$, \ then one can define
\ $ \delta$ \ standing on \ $\alpha _{0}^{\downarrow }$\
 and can see that all conditions of this lemma are
fulfilled and therefore it implies the contradiction:

\quad \quad \quad \quad \quad \quad \ \ $Od ( S_{\tau ^{\ast }} )
\leq \delta <OT
    ( S_{\tau_{0}} ) \leq Od ( S_{\tau ^{\ast }} )$.

 \noindent
Hence, \ $S_f$ is nonmonotone and it constitutes the required
contradiction.

 \noindent
The case \ $\delta =\chi ^{\ast +}$ \ can be eliminated by certain
transformation of matrix functions \ $S_{f}^{<\alpha_1}$. \

\end{document}